\documentclass[11pt]{article} \usepackage{amsfonts} \usepackage{bbm} \usepackage{amscd} \usepackage{xypic} \usepackage{amsmath} \usepackage{amssymb} \usepackage{amsthm} \usepackage{bbm} \usepackage{CJK} \usepackage{fancyhdr} \usepackage{hyperref} \usepackage{indentfirst} \usepackage{latexsym} \usepackage{mathrsfs}
 \allowdisplaybreaks 
\usepackage{color}

\usepackage[top=1in,bottom=1in,left=1.25in,right=1.25in]{geometry} \textwidth 15.5cm \textheight 22cm \oddsidemargin 0in \def\<{\langle} \def\>{\rangle}
\def\a{\alpha} \def\b{\beta} \def\ci{\circ} \def\c{\cdot}      \def\i{\iota}
   \def\r{\rho} \def\lr{\longrightarrow}  
\def\o{\otimes}       \def\vp{\varphi} 
\def\<{\langle} \def\>{\rangle}   \date{} \begin{document} \renewcommand{\baselinestretch}{1.2} \renewcommand{\arraystretch}{1.0}
\title{\bf New Braided $T$-Categories over Hopf (co)quasigroups} \date{} \author{{\bf Wei Wang and
        Shuanhong Wang\footnote {Corresponding author:Shuanhong Wang, E-mail: shuanhwang2002@yahoo.com}
}\\ {\small Department of Mathematics, Southeast University}\\ {\small Jiangsu Nanjing 210096, P. R. CHINA}\\}
 \maketitle
\begin{center} \begin{minipage}{14.cm} \begin{center}{\bf ABSTRACT}\end{center}

Let $H$ be a Hopf quasigroup with bijective antipode and let $Aut_{HQG}(H)$ be the set of all Hopf quasigroup automorphisms of $H$. We introduce a category
${_{H}\mathcal{YDQ}^{H}}(\alpha,\beta)$ with $\alpha,\beta\in Aut_{HQG}(H)$ and construct a braided $T$-category $\mathcal{YDQ}(H)$ having all the
categories ${_{H}\mathcal{YDQ}^{H}}(\alpha,\beta)$ as components.

 \vskip 0.5cm

{\bf Key words}:  Hopf quasigroup; braided $T$-category; Quasi $(\alpha,\beta)$-Yetter-Drinfeld category.
 \vskip 0.5cm
 {\bf Mathematics Subject Classification 2010:} 16T05.
\end{minipage} \end{center}

\section*{ Introduction} \def\theequation{0. \arabic{equation}} \setcounter{equation} {0} \hskip\parindent

To understand the structure and relevant properties of the algebraic 7-sphere, Klim and Majid in \cite{KM10} introduced the notion of Hopf quasigroups.
Which is the ¡®linearise¡¯ of quasigroups, here Hopf quasigroups are not associative but the lack of this property is compensated by some axioms involving
the antipode $S$. The concept of Hopf quasigroup is a particular case of the notion of unital coassociative $H$-bialgebra introduced in \cite{P07}, and
Hopf quasigroup includes the example of an enveloping algebra $U(L)$ of a Mal¡¯tsev algebra (see \cite{KM10}) as well as the notion of quasigroup algebra
$RL$ of an I.P. loop $L$. In particular, Hopf quasigroups unify I.P. loops and Mal¡¯tsev algebras as the same as Hopf algebras unified groups and Lie
algebras. \\

Turaev in \cite{V08, V94} generalized quantum invariants for 3-manifolds to the case of a 3-manifold $M$ endowed with a homotopy class of maps
$M\rightarrow K(G,1)$, where $G$ is a group. And braided $T$-categories which are braided monoidal categories in Freyd-Yetter categories of crossed
$G$-sets(see in \cite{PD89}) play a key role of constructing these homotopy invariants.\\

Based on these two notions and structures, the aim of this paper is to construct classes of new braided $T$-categories over Hopf quasigroups. We first
introduce the concept of an $(\alpha,\beta)$-Yetter-Drinfeld quasimodule over a Hopf quasigroup, which is a generalization of a Yetter-Drinfeld quasimodule
over a Hopf quasigroup (see \cite{AFGS15}) , then we construct new examples of braided $T$-categories, which generalize the construction of braided
$T$-cateogries over Hopf algebras given by Panaite and Staic (see in \cite{PM07}).\\

Let $H$ be a Hopf quasigroup with bijective antipode. We denote by $Aut_{HQG}(H)$ the set of all Hopf quasigroup automorphisms $\alpha$ of $H$ that
satisfying $S\circ \alpha= \alpha \circ S$, and we consider $G$ a certain crossed product group $Aut_{HQG}(H)\times Aut_{HQG}(H)$.\\

In Section 1, we recall definitions and basic results related to Hopf quasi-groups, Yetter-Drinfeld quasimodules over Hopf quasigroups and braided
$T$-categories. In Section 2, we introduce a class of new categories $ _{H}\mathcal{YDQ}^{H}(\alpha, \beta)$ (see Definition 2.1) of
$(\alpha,\beta)$-Yetter-Drinfeld quasimodules associated with $\alpha,\beta\in Aut_{HQG}(H)$. Then in Section 3,  we prove $\mathcal{YDQ}(H)$ is a monoidal
cateogory and then construct a class of new braided $T$-categories $\mathcal{YDQ}(H)$ in the sense of Turaev\cite{V08}.\\

\section*{1. Preliminaries} \def\theequation{1. \arabic{equation}} \setcounter{equation} {0} \hskip\parindent

Throughout, let $k$ be a fixed field. Everything is over $k$ unless
 otherwise specified.  We
 refer the readers to the book of Sweedler \cite{M69}
  for the relevant concepts on the general theory of Hopf
 algebras.  Let $(C, \Delta )$ be a coalgebra. We use the simplified Sweedler-Heyneman's notation for
 $\Delta $ as follows:
 \begin{eqnarray*}
 \Delta (c)=c_1\otimes c_2,
 \end{eqnarray*}
 for all $c\in C$.

 \vskip 0.5cm
 {\bf 1.1. Hopf (co)quasigroups.}
 \vskip 0.5cm

 We first recall that an (inverse property) quasigroup is a Set $Q$ with a product,
 an identity $e$ and for each $s\in Q$, there is an element $s^{-1}\in Q$ such that
 \begin{eqnarray*}
 s^{-1}(st)=t,\quad (ts)s^{-1}=t, \quad \forall t\in Q.
 \end{eqnarray*}
 A quasigroup is {\sl Moufang} if $s(t(sr))=((st)s)r$ for all $s,t,r\in Q$.
 In \cite{KM10}, Klim and Majid linearised these notions to Hopf quasigroups in the same way that Hopf algebras linearises the notions of groups.

 A {\sl Hopf quasigroup} is a unital algebra $H$ (possibly-nonassociative), equipped with algebra maps $\Delta:H\rightarrow H\otimes H$ and
 $\varepsilon: H\rightarrow k$ forming a coassociative coalgebra and a map $S:H\rightarrow H$ such that
 \begin{eqnarray*}
 && S(h_{1})(h_{2}g)=g\varepsilon(h)= h_{1}(S(h_{2})g),\\
 && (gS(h_{1}))h_{2}=g\varepsilon(h)=(gh_{1})S(h_{2}),
 \end{eqnarray*}
 for all $h,g\in H$. In this notation the Hopf quasigroup $H$ is called {\sl Moufang} if
 \begin{eqnarray*}
  h_{1}(g(h_{2}f))= ((h_{1}g)h_{2})f,\quad \forall h,g,f\in H.
 \end{eqnarray*}
 And a Hopf quasigroup $H$ is called {\sl flexible} if
 \begin{eqnarray*}
   h_{1}(gh_{2})=(h_{1}g)h_{2},\quad \forall g,h\in H.
 \end{eqnarray*}

 We know here the conditions of antipode $S$ are stronger than the usual Hopf algebra antipode axioms and then compensate for $H$ nonassociative.
 For instance, $S$ is antimultiplicative and anticomultiplicative in the sense
 \begin{eqnarray*}
 S(hg)=S(g)S(h),\quad \Delta(S(h))=S(h_{2})\otimes S(h_{1}),
 \end{eqnarray*}
 for all $g,h\in H$. If we linearise an (inverse property) quasigroup $Q$ to a Hopf quasigroup algebra $kQ$ with grouplike coproduct on elements of $Q$
 and
 linear extension of the product and inverse, and the Hopf quasigroup is {\sl Moufang} if $Q$ is {\sl Moufang}.

 A (left) {\sl Hopf quasimodule} (see \cite{BJ12}) over a Hopf quasigroup $H$ is a vector space $M$ equipped with a structure $\cdot: H\otimes M\rightarrow
 M$ such that
 \begin{eqnarray*}
 1_{H}\cdot m&=&m,\\
 h_{1}\cdot(S(h_{2})\cdot m)&=&S(h_{1})\cdot (h_{2}\cdot m)\,\,=\,\,\varepsilon(h)m.
 \end{eqnarray*}
 for all $h\in H$ and $m\in M$.

 Dually we also know the notion of the {\sl Hopf coquasigroup}, that $H$ is an associative algebra and a couintal coalgebra (the comultiplication
 is possibly noncoassociative), and a map $S:H\rightarrow H$ such that
 \begin{eqnarray*}
 S(h_{1})h_{21}\otimes h_{22}&=&1_{H}\otimes h\,\,=\,\,h_{1}S(h_{21})\otimes h_{22},\\
 h_{11}\otimes S(h_{12})h_{2}&=&h\otimes 1_{H}\,\,=\,\,h_{11}\otimes h_{12}S(h_{2}),
 \end{eqnarray*}
 for all $h\in H$. Similarly we call a Hopf coquasigroup {\sl Flexible} if $h_{1}h_{22}\otimes h_{21}=h_{11}h_{2}\otimes h_{12}, \forall h\in H$,
 and {\sl Moufang} if $h_{1}h_{221}\otimes h_{21}\otimes h_{222}=h_{111}h_{12}\otimes h_{112}\otimes h_{2}, \forall h\in H$.
 \\

 \vskip 0.5cm
 {\bf 1.2. Yetter-Drinfeld quasimodules over a Hopf quasigroup.}
 \vskip 0.5cm

 Let $H$ be a Hopf quasigroup, in \cite{AFGS15}, authors gave the notion of left-left Yetter-Drinfeld quasimodule over $H$.
 Similarly, we say that $M=(M,\cdot, \rho)$ is a left-right Yetter-Drinfeld quasimodule over $H$ if
 $(M,\cdot)$ is a left $H$-quasimodule and $(M,\rho)$ is a right $H$-comodule which satifies the following:
 \begin{eqnarray}
 (h_{2}\cdot m)_{(0)}\otimes (h_{2}\cdot m)_{(1)}h_{1}&=&h_{1}\cdot m_{(0)}\otimes h_{2}m_{(1)},\\
 m_{(0)}\otimes m_{(1)}(hg)&=&m_{(0)}\otimes (m_{(1)}h)g,\\
 m_{(0)}\otimes h(m_{(1)}g)&=&m_{(0)}\otimes (hm_{(1)})g,
 \end{eqnarray}
 for all $h,g\in H$ and $m\in M$. Here we use the notation of right $H$-comodule by $\rho(m)=m_{(0)}\otimes m_{(1)}, \forall m\in M$.

 Let $M$ and $N$ be two left-right Yetter-Drinfeld quasimodules over $H$. We call morphism $f:M\rightarrow N$ a left-right Yetter-Drinfeld quasimodule
 morphism if
 $f$ is both a left $H$-quasimodule morphism and a right $H$-comodule morphism. We use $_{H}\mathcal{YDQ}^{H}$ denote the category of left-right
 Yetter-Drinfeld quasimodules
 over $H$. Moreover, if we assume $H$ is associative, that $H$ is a Hopf algebra, then conditions (1.2) and (1.3) become trivial. In this case, if $M$ is a
 left $H$-module we
 obtain the classical definition of left-right Yetter-Drinfeld module over a Hopf algebra.

 \vskip 0.5cm
 {\bf 1.3. Braided $T$-categories.}
 \vskip 0.5cm

A {\sl monoidal category} $\mathcal{C}=(\mathcal{C},\mathbb{I},\otimes,a,l,r)$
 is a category $\mathcal{C}$ endowed with a functor
 $\otimes: \mathcal{C}\times\mathcal{C}\rightarrow\mathcal{C}$
 (the {\sl tensor product}), an object $\mathbb{I}\in \mathcal{C}$
 (the {\sl tensor unit}), and natural isomorphisms $a$
 (the {\sl associativity constraint}),$l$ (the {\sl left unit constraint})
 and $r$ (the {\sl right unit constraint}),
 such that for all $U,V,W,X\in \mathcal{C},$
 the {\sl associativity pentagon}
 $a_{U,V,W\otimes X}\circ a_{U\otimes V,W,X}
 =(U\otimes a_{V,W,X})\circ a_{U,V\otimes W,X}\circ
 (a_{U,V,W}\otimes X)$ and
 $(U\otimes l_V)\circ(r_U\otimes V)=a_{U,I,V}$ are satisfied.
 A monoidal categoey $\mathcal{C}$ is {\sl strict} when all
  the constraints are identities.
\\

  Let $G$ be a group and let $Aut(\mathcal{C})$ be the group of
  invertible strict tensor functors from $\mathcal{C}$ to itself.
  A category $\mathcal{C}$ over $G$
  is called a {\sl crossed category } if it satisfies the following:

 $\bullet$  $\mathcal{C}$ is a monoidal category;

 $\bullet$  $\mathcal{C}$ is disjoint union of a family
 of subcategories $\{\mathcal{C}_{\alpha }\}_{\alpha \in G}$, and for any $U\in \mathcal{C}_{\alpha }$,
 $V\in \mathcal{C}_{\beta }$, $U\otimes V\in \mathcal{C}_{\a \beta }$.
 The subcategory $\mathcal{C}_{\alpha}$ is called the $\alpha$th component of $\mathcal{C}$;

 $\bullet$  Consider a group homomorphism $\varphi : G\rightarrow Aut(\mathcal{C})$, $\b \mapsto \varphi _{\b }$, and
  assume that$\varphi_{\beta }(\varphi_{\alpha })=\varphi_{\beta\alpha\beta^{-1}}$,for all $\alpha,\beta\in G$.
  The functors $\varphi_{\beta }$  are called conjugation isomorphisms.

 Furthermore, $\mathcal{C}$ is called strict when it is
 strict as a monoidal category.
\\

 {\sl Left index notation}: Given $\a \in G$
 and an object $V\in \mathcal{C}_{\a }$, the functor $\vp _{\a }$
 will be denoted by ${}^V( \cdot )$, as in Turaev \cite{V08} or
 Zunino \cite{M04}, or even ${}^{\a }( \cdot )$.
 We use the notation ${}^{\overline{V}}( \cdot )$
 for ${}^{\a ^{-1}}( \cdot )$. Then we have
 ${}^V id_U=id_{{} V^U}$ and
 ${}^V(g\circ f)={}^Vg\circ {}^Vf$.
 Since the conjugation $\vp : G\lr Aut(\mathcal{C})$ is a
 group homomorphism, for all $V, W\in \mathcal{C}$, we have ${}^{V\o W}( \cdot )
 ={}^V({}^W( \cdot ))$ and ${}^\mathbb{I}( \cdot )={}^V({}^{\overline{V}}( \cdot ))
 ={}^{\overline{V}}({}^V( \cdot ))=id_\mathcal{C}$. Since, for all
 $V\in \mathcal{C}$, the functor ${}^V( \cdot )$ is strict, we have
 ${}^V(f\o g)={}^Vf\o {}^Vg$, for any morphisms $f$ and $g$ in $\mathcal{C}$,
  and ${}^V\mathbb{I}=\mathbb{I}$.
\\

 A {\sl braiding} of a crossed category $\mathcal{C}$ is
 a family of isomorphisms $({c=c_{U,V}})_{U,V}\in \mathcal{C}$,
 where $c_{U,V}: U\otimes V\rightarrow {}^UV\otimes U$
 satisfying the following conditions:\\
 (i)~For any arrow $f\in \mathcal{C}_{\a }(U, U')$ and
 $g\in \mathcal{C}(V, V')$,
 $$
 (({}^{\a }g)\o f)\circ c _{U, V}=c _{U' V'}\circ (f\o g).
 $$
 (ii)~For all $ U, V, W\in \mathcal{C},$ we have
 $$
 c _{U\o V, W}=a_{{}^{U\o V}W, U, V}\circ (c _{U, {}^VW}\o
 id_V)\circ a^{-1}_{U, {}^VW, V}\circ (\i _U\o c _{V, W})
  \circ a_{U, V, W},
  $$
    $$c _{U, V\o W}=a^{-1}_{{}^UV, {}^UW, U}
 \circ (\i _{({}^UV)}\o c _{U, W})\circ a_{{}^UV, U, W}\ci
 (c _{U, V}\o \i_W)\circ a^{-1}_{U, V, W},$$
 where $a$ is the natural isomorphisms in the tensor category
 $\mathcal{C}$.\\
 (iii)~For all $ U, V\in \mathcal{C}$ and $\b\in G$,
 $$ \vp _{\b }(c_{U, V})=c_{\vp _{\b }(U), \vp_{\b }(V)}.
 $$

 A crossed category endowed with a braiding is called
 a {\sl braided $T$-category}.\\

\section*{2. $(\alpha, \beta)$- Yetter-Drinfeld quasimodules over a Hopf quasigroup} \def\theequation{2. \arabic{equation}} \setcounter{equation} {0}
\hskip\parindent

 In this section, we will define the notion of a Yetter-Drinfeld quasimodule over
 a Hopf quasigroup that is twisted by two Hopf quasigroup automorphisms as well as the notion of a Hopf quasi-entwining structure
 and how to obtain such structure from automorphisms of Hopf quasigroups.
\\

In what follows, let  $H$ be a Hopf quasigroup with the bijective antipode $S$  and let ${\sl Aut}_{HQG}(H)$ denote the set of all automorphisms of a Hopf
quasigroup $H$. \\

{\bf Definition 2.1.} Let $\alpha, \beta\in  {\sl Aut}_{HQG}(H)$. A left-right {\sl $(\alpha, \beta)$-Yetter-Drinfeld quasimodule} over $H$ is a vector
space $M$, such that $M$ is a left $H$-quasimodule(with notation $h\otimes m\mapsto h\cdot m$) and a right $H$-comodule(with notation $M\rightarrow
M\otimes H$, $m\mapsto m_{(0)}\otimes m_{(1)}$) and with the following compatibility condition:
 \begin{equation}
\r (h\cdot m)=h_{21}\cdot m_{(0)}\otimes (\beta(h_{22})m_{(1)})\alpha(S^{-1}(h_{1})), \end{equation} for all $h\in H$ and $m\in M$. We denote by $
_{H}\mathcal{YDQ}^{H}(\alpha, \beta)$ the category of left-right $(\alpha, \beta)$-Yetter-Drinfeld quasimodules, morphisms being both $H$-linear and
$H$-colinear maps. \\

 {\bf Remark.} Note that, $\alpha$ and $\beta$ are bijective,
 algebra morphisms, coalgebra morphisms and commute with $S$.
 \\

{\bf Proposition 2.2.} One has that Eq.(2.1) is equivalent to
 the following equations:
\begin{eqnarray} h_{1}\c m_{(0)}\o \beta(h_{2})m_{(1)}=(h_{2}\cdot m)_{(0)}\otimes (h_{2}\cdot m)_{(1)}\alpha(h_{1}). \end{eqnarray}

{\bf Proof.}  To prove this propostion, we need to use the property of antipode of a Hopf quasigroup. that is,
$(gS(h_{1}))h_{2}=(gh_{1})S(h_{2})=g\varepsilon(h)$, for all $g,h\in H$.

Eq.(2.1)$\Longrightarrow$ Eq.(2.2).

We first do calculation as follows: \begin{eqnarray*}
 &&(h_{2}\cdot m)_{(0)}\otimes (h_{2}\cdot m)_{(1)}\alpha(h_{1})\\
 &\stackrel{(2.1)}{=}& h_{22}\cdot m_{(0)}\otimes ((\beta(h_{23})m_{(1)})\alpha(S^{-1}(h_{21})))\alpha(h_{1})\\
 &=& h_{22}\cdot m_{(0)}\otimes ((\beta(h_{23})m_{(1)})\alpha(S^{-1}(h_{21})))\alpha(h_{1})\\
 &=& h_{2}\cdot m_{(0)}\otimes ((\beta(h_{3})m_{(1)})\varepsilon(h_{1})\\
 &=& h_{1}\cdot m_{(0)}\otimes \beta(h_{2})m_{(1)}.
\end{eqnarray*}

 For Eq.(2.2) $\Longrightarrow$ Eq.(2.1), we have
\begin{eqnarray*} &&h_{2}\cdot m_{(0)}\otimes (\beta(a_{3})m_{(1)})\alpha(S^{-1}(a_{1}))\\ &\stackrel{(2.2)}{=}& (h_{3}\cdot m)_{(0)}\otimes ((a_{3}\cdot
m)_{(1)}\alpha(a_{2}))\alpha(S^{-1}(a_{1}))\\ &=& (a_{2}\cdot m)_{(0)}\otimes (a_{2}\cdot m)_{(1)}\varepsilon(a_{1})\\ &=& (a\cdot m)_{(0)}\otimes (a\cdot
m)_{(1)}. \end{eqnarray*} This finishes the proof. \hfill $\blacksquare$ \\

{\bf Example 2.3.}  For $\beta\in {\sl Aut}_{HQG}(H)$, define $H_{\beta}=H$ as vector space over a field $k$, with regular right $H$-comdule structure and
left $H$-quasimodule structure given by $h\cdot h'=(\beta(h_{2})h')S^{-1}(h_{1})$, for all $h,h'\in H$. And more generally, if $\alpha,\beta\in {\sl
Aut}_{HQG}(H)$, define $H_{\alpha,\beta}=H$ with regular right $H$-comodule structure and left $H$-module structure given by $h\cdot
h'=(\beta(h_{2})h')\alpha(S^{-1}(h_{1}))$, for all $h,h'\in H$. If $H$ is flexible, then we get $H_{\beta}\in  _{H}\mathcal{YDQ}^{H}(id, \beta)$ and if we
add $H$ is $(\alpha,\beta)$-flexible, that is \begin{eqnarray*} \alpha(h_{1})(g\beta(h_{2}))=(\alpha(h_{1})g)\beta(h_{2}), \end{eqnarray*} for all $g,h\in
H$and $\alpha,\beta\in {\sl Aut}_{HQG}(H)$ . Then $H_{\alpha,\beta}\in  _{H}\mathcal{YDQ}^{H}(\alpha, \beta)$. \\

 Let $\alpha,\beta\in {\sl Aut}_{HQG}(H)$. As defined in \cite{PM07}, an $H$-bicomodule algebra $H(\alpha,\beta)$ as follows; $H(\alpha,\beta)=H$ as
 algebras, with comodule structures
 \begin{eqnarray*}
&& H(\alpha,\beta)\rightarrow H\otimes H(\alpha,\beta),\quad h\mapsto h_{[-1]}\otimes h_{[0]}=\alpha(h_{1})\otimes h_{2},\\ &&  H(\alpha,\beta)\rightarrow
H(\alpha,\beta)\otimes H,\quad h\mapsto h_{<0>}\otimes h_{<1>}=h_{1}\otimes \beta(h_{2}).
 \end{eqnarray*}
 Then we also consider the  Yetter-Drinfeld quasimodules like $_{H(\alpha,\beta)}\mathcal{YDQ}^{H}(H)$.

{\bf Propostion 2.4.} $_{H}\mathcal{YDQ}^{H}(\alpha, \beta)=_{H(\alpha,\beta)}\mathcal{YDQ}^{H}(H)$.

{\bf Proof.} Left to readers. \\

\section*{3. A BRAIDED $T$-CATEGORY $\mathcal {YDQ}(H)$} \def\theequation{3. \arabic{equation}} \setcounter{equation} {0} \hskip\parindent

In this section, we will construct a class of new braided $T$-categories
 $\mathcal {YDQ}(H)$ over any Hopf quasigroup $H$ with bijective antipode.Here for $\alpha \in {\sl Aut}_{HQG}(H)$,
 the ${_{H}}\mathcal {YDQ}^{H}(\alpha,\beta)$ is the object of $\mathcal {YDQ}(H).$
 \\

 Let $M\in {_{H}}\mathcal {YDQ}^{H}(\alpha,\beta)$
 ,$N\in {_{H}}\mathcal {YDQ}^{H}(\gamma,\delta)$, with
 $\alpha,\beta,\gamma,\delta \in {\sl Aut}_{HQG}(H)$.
\\

{\bf Proposition 3.1.} If $M\in {_{H}}\mathcal {YDQ}^{H}(\alpha,\beta)$
 ,$N\in {_{H}}\mathcal {YDQ}^{H}(\gamma,\delta)$, with
 $\alpha,\beta,\gamma,\delta \in {\sl Aut}_{HQG}(H)$., then $M \otimes N
 \in {_{H}}\mathcal {YDQ}^{H}(\alpha\gamma, \delta\gamma^{-1}\beta\gamma)$ with structures as follows:
\begin{eqnarray*} &&h\c (m \otimes n)=\gamma (h_{1})\c m \otimes \gamma^{-1}\beta\gamma(h_{2})\c n,\\ &&m\otimes n \mapsto (m_{(0)}\otimes n_{(0)})\otimes
n_{(1)}m_{(1)}. \end{eqnarray*} for all $m\in M,n\in N$ and $h\in H.$ \\

{\bf Proof.} Let $h,g\in H$ and $m\otimes n\in M\otimes N$. We can prove $1\cdot (m\otimes n)=m\otimes n$ and $h_{1}\cdot(S(h_{2})\cdot (m\otimes
n))=S(h_{1})\cdot(h_{2}\cdot (m\otimes n))=(m\otimes n)$, straightforwardly.

This shows that $M \otimes N$ is a left $H$-quasimodule, the right $H$-comodule condition is straightforward to check.

Next, we compute the compatibility condition as follows: \begin{eqnarray*} && (h_{2}\cdot (m\otimes n))_{(0)}\otimes (h_{2}\cdot (m\otimes
n))_{(1)}\alpha\gamma(h_{1})\\ &=& (\gamma(h_{2})\cdot m\otimes \gamma^{-1}\beta\gamma(h_{3})\cdot n)_{(0)}\otimes (\gamma(h_{2})\cdot m\otimes
\gamma^{-1}\beta\gamma(h_{3})\cdot n)_{(1)}\alpha\gamma(h_{1})\\ &=& (\gamma(h_{2})\cdot m)_{(0)}\otimes (\gamma^{-1}\beta\gamma(h_{3})\cdot
n)_{(0)}\otimes (\gamma^{-1}\beta\gamma(h_{3})\cdot n)_{(1)}((\gamma(h_{2})\cdot m)_{(1)}\alpha\gamma(h_{1}))\\ &=& \gamma(h_{1})\cdot m_{(0)}\otimes
(\gamma^{-1}\beta\gamma(h_{3})\cdot n)_{(0)}\otimes ((\gamma^{-1}\beta\gamma(h_{3})\cdot n)_{(1)}\gamma\gamma^{-1}\beta\gamma(h_{2}))m_{(1)}\\ &=&
\gamma(h_{1})\cdot m_{(0)}\otimes \gamma^{-1}\beta\gamma(h_{2})\cdot n_{(0)}\otimes \delta\gamma^{-1}\beta\gamma(h_{3})(n_{(1)}m_{(1)})\\ &=&
h\dot(m_{(0)}\otimes n_{(0)})\otimes \delta\gamma^{-1}\beta\gamma(h_{2})(n_{(1)}m_{(1)}). \end{eqnarray*} Thus $M \otimes N \in {_{H}}\mathcal
{YDQ}^{H}(\alpha\gamma, \delta\gamma^{-1}\beta\gamma)$.
 \hfill $\blacksquare$
\\

{\bf Remark. } Note that, if $M \in {_{H}}\mathcal {YDQ}^{H}(\alpha,\beta),\,\, N \in {_{H}}\mathcal {YDQ}^{H}(\gamma,\delta)$ and $P\in {_{H}}\mathcal
{YDQ}^{H}(\mu, \nu)$, then $(M\otimes N)\otimes P = M\otimes (N\otimes P)$ as objects in ${_{H}}\mathcal
{YDQ}^{H}(\alpha\gamma\mu,\nu\mu^{-1}\delta\gamma^{-1}\beta\gamma\mu).$ \\

 Denote $G={\sl Aut}_{HQG}(H)
 \times {\sl Aut}_{HQG}(H)$
 a group with multiplication as follows:
  for all $\alpha,\beta, \gamma, \delta \in {\sl Aut}_{HQG}(H)$,
\begin{equation} (\alpha,\beta)\ast (\gamma, \delta)=(\alpha\gamma, \delta\gamma^{-1}\beta\gamma). \end{equation}
 The unit of this group is $(id,id)$ and $(\alpha,\beta)^{-1}=(\alpha^{-1},
 \alpha\beta^{-1}\alpha^{-1})$.\\

The above proposition means that if $M \in {_{H}}\mathcal {YDQ}^{H}(\alpha,\beta)$
 and $N\in {_{H}}\mathcal {YDQ}^{H}(\gamma, \delta)$, then
 $M \otimes N \in {_{H}}\mathcal {YDQ}^{H}((\alpha,\beta)\ast (\gamma, \delta)).$
\\

{\bf Proposition 3.2.} Let $N \in {_{H}}\mathcal {YDQ}^{H}(\gamma, \delta)$
 and $(\alpha,\beta)\in G$. Define ${}^{(\alpha,\beta)}N=N$ as
 vector space, with structures: for all $n\in N$ and $h\in H.$
$$h\rhd n=\gamma^{-1}\beta\gamma\alpha^{-1}(h)\c n,$$ \begin{equation} n\mapsto n_{<0>}\otimes n_{<1>}=n_{(0)}\otimes \alpha\beta^{-1}(n_{(1)}).
\end{equation} Then $${}^{(\alpha,\beta)}N \in {_{H}}\mathcal {YDQ}^{H}((\alpha,\beta)\ast (\gamma,\delta)\ast (\alpha,\beta)^{-1}),$$ where
$((\alpha,\beta)\ast (\gamma,\delta)\ast (\alpha,\beta)^{-1})=(\alpha\gamma\alpha^{-1}, \alpha\beta^{-1}\delta\gamma^{-1}\beta\gamma\alpha^{-1})$ as an
element in $G$.

{\bf Proof.} Obviously, the equations above define a  quasi-module and a comodule action of $N$. In what follows, we show the compatibility condition:
\begin{eqnarray*} &&(h\rhd n)_{<0>}\otimes(h\rhd n)_{<1>}\\ &=& (\gamma^{-1}\beta\gamma\alpha^{-1}(h)\c n)_{(0)}\otimes
\alpha\beta^{-1}((\gamma^{-1}\beta\gamma\alpha^{-1}(h)\c n)_{(1)})\\ &=& \gamma^{-1}\beta\gamma\alpha^{-1}(h)_{2}\cdot n_{(0)}\otimes
\alpha\beta^{-1}((\delta(\gamma^{-1}\beta\gamma\alpha^{-1}(h)_{3})n_{(1)}) \gamma(S^{-1}(\gamma^{-1}\beta\gamma\alpha^{-1}(h)_{1})))\\ &=& h_{2}\rhd
n_{(0)}\otimes (\alpha\beta^{-1}\delta\gamma^{-1}\beta\gamma\alpha^{-1}(h_{3})\alpha\beta^{-1}(n_{(1)}))\alpha\gamma\alpha^{-1}(S^{-1}(h_{1}))\\ &=&
h_{2}\rhd n_{(0)}\otimes (\alpha\beta^{-1}\delta\gamma^{-1}\beta\gamma\alpha^{-1}(h_{3})n_{<1>})\alpha\gamma\alpha^{-1}(S^{-1}(h_{1})) \end{eqnarray*} for
all $n\in N$ and $h\in H,$ that is ${}^{(\alpha,\beta)}N \in {_{H}}\mathcal
{YDQ}^{H}(\alpha\gamma\alpha^{-1},\alpha\beta^{-1}\delta\gamma^{-1}\beta\gamma\alpha^{-1})$ \hfill $\blacksquare$ \\

{\bf Remark.} Let $M \in {_{H}}\mathcal {YDQ}^{H}(\alpha, \beta),
 \ \ N \in {_{H}}\mathcal
{YDQ}^{H}(\gamma,\delta),\,\, \mbox {and}\, (\mu,\nu)\in G$. Then by the
 above proposition, we have:
 $$
 {}^{(\alpha, \beta)\ast (\mu,\nu)}N={}^{(\alpha, \beta)}({}^{(\mu,\nu)}N),
$$ as objects in $_{H}\mathcal {YDQ}^{H}(\alpha \mu\gamma \mu^{-1}\alpha^{-1}, \alpha\beta^{-1}\mu\nu^{-1}\delta\gamma^{-1}\nu\mu^{-1}\beta \mu\gamma
\mu^{-1}\alpha^{-1})$ and $$ {}^{(\mu,\nu)}(M\otimes N)= {}^{(\mu,\nu)}M \otimes {}^{(\mu,\nu)}N, $$ as objects in $_{H}\mathcal {YDQ}^{H}(\mu\alpha \gamma
\mu^{-1},
 \mu\nu^{-1}\delta\gamma^{-1}\beta\alpha^{-1}\nu\alpha \gamma \mu^{-1})$.
\\

 {\bf Proposition 3.3.} Let $M \in {_{H}}\mathcal {YDQ}^{H}(\alpha,\beta)$,
 $N \in {_{H}}\mathcal {YDQ}^{H}(\gamma,\delta)$ and $ P \in {_{H}}\mathcal {YDQ}^{H}(\mu,\nu)$, take
  ${}^{M}N={}^{(\alpha,\beta)}N$ as explained in Subsection 1.3.
  Define a map $c_{M, N}: M \otimes N \rightarrow {}^{M}N \otimes M$ by
\begin{equation}
 c_{M,N}(m\otimes
 n)=n_{(0)}\otimes \beta^{-1}(n_{(1)})\cdot m.
\end{equation} for all $m\in M,n\in N.$ Then $c_{M, N}$ is both an $H$-module map and an $H$-comodule map,
 and satisfies the following formulae:
\begin{equation} c_{M\otimes N,P}=(c_{M,^{N}P}\otimes id_{N})\circ (id_{M}\otimes c_{N,P}), \end{equation} \begin{equation}
 c_{M,N\otimes P}=(id_{^{M}N}\otimes c_{M,P})\circ(c_{M,N}\otimes id_{P}).
\end{equation} \\

{\bf Proof.} First, we prove that $c_{M,N}$ is an $H$-module map.
 Take  $h\cdot(m\otimes n)=\gamma(h_{1})\c m\otimes \gamma^{-1}\beta\gamma(h_{2})\c n$
 as explained in Proposition 3.1.
 \begin{eqnarray*}
&&c_{M,N}(h\cdot(m\otimes n))\\ &=& c_{M,N}(\gamma(h_1)\cdot m\otimes \gamma^{-1}\beta\gamma(h_2)\cdot n)\\ &=& (\gamma^{-1}\beta\gamma(h_{2})\cdot
n)_{(0)}\otimes \beta^{-1}((\gamma^{-1}\beta\gamma(h_{2})\cdot n)_{(1)})\cdot (\gamma(h_{1})\cdot m)\\ &=& (\gamma^{-1}\beta\gamma(h_{2})\cdot
n)_{(0)}\otimes (\beta^{-1}((\gamma^{-1}\beta\gamma(h_{2})\cdot n)_{(1)})\gamma(h_{1}))\cdot m\\ &=& \gamma^{-1}\beta\gamma(h_{3})\cdot n_{(0)}\otimes
((\beta^{-1}\delta\gamma^{-1}\beta\gamma(h_{4})\beta^{-1}(n_{(1)}))S^{-1}(\gamma(h_{2})))\gamma(h_{1})\cdot m\\ &=& \gamma^{-1}\beta\gamma(h_{2})\cdot
n_{(0)}\otimes (\beta^{-1}\delta\gamma^{-1}\beta\gamma(h_{3})\beta^{-1}(n_{(1)}))\varepsilon(h_{1})\cdot m\\ &=& \gamma^{-1}\beta\gamma(h_{1})\cdot
n_{(0)}\otimes (\beta^{-1}\delta\gamma^{-1}\beta\gamma(h_{2})\beta^{-1}(n_{(1)}))\cdot m, \end{eqnarray*} on the other side, we have \begin{eqnarray*}
  h\cdot c_{M,N}(m\otimes n) &=& h\cdot (n_{(0)}\otimes \beta^{-1}(n_{(1)})\cdot m)\\
  &=& \alpha(h_{1})\rhd n_{(0)} \otimes (\beta^{-1}\delta\gamma^{-1}\beta\gamma(h_{2})\beta^{-1}(n_{(1)}))\cdot m,
\end{eqnarray*} similarly we can check that $c_{M, N}$ is an $H$-comodule map.

Finally we will check Eqs.(3.4) and (3.5). Using equations $^M(^N P)=^{M\otimes N}P$ and $^M(N\otimes P)=^M N\otimes ^M P$ we have \begin{eqnarray*}
&&(c_{M,^{N}P}\otimes id_{N})\circ (id_{M}\otimes c_{N,P})(m\otimes n\otimes p)\\ &=& c_{M,^{N}P}(m\otimes p_{(0)})\otimes \delta^{-1}(p_{(1)})\cdot n\\
&=& p_{(0)(0)}\otimes \beta^{-1}\gamma\delta^{-1}(p_{(0)(1)})\cdot m\otimes \delta^{-1}(p_{(1)})\cdot n\\ &=& p_{(0)}\otimes
\beta^{-1}\gamma\delta^{-1}(p_{(1)(1)})\cdot m\otimes \delta^{-1}(p_{(1)(2)})\cdot n\\ &=& p_{(0)}\otimes
\gamma^{-1}\beta^{-1}\gamma\delta^{-1}(p_{(1)})\cdot (m\otimes n)\\ &=& c_{M\otimes N,P}(m\otimes n\otimes p). \end{eqnarray*} Similar we can check the
equation 3.5, that ends the proof.
 \hfill $\blacksquare$
\\

{\bf Lemma 3.5.} The map $c_{M,N}$ defined by $c_{M,N}(m\otimes n)=n_{(0)} \otimes \beta^{-1}(n_{(1)})\cdot m$ is bijective; with inverse
 $${c}_{M,N}^{-1}(n\otimes m)=\beta^{-1}(S(n_{(1)}))\cdot m \otimes n_{(0)}.$$

{\bf Proof.} First, we prove $c_{M,N}\circ c_{M,N}^{-1}=id$. For all $m\in M, n\in
  N$, we have
 \begin{eqnarray*}
&&c_{M,N}\circ c_{M,N}^{-1}(n \otimes m)\\ &=&c_{M,N}(\beta^{-1}(S(n_{(1)}))\cdot m\otimes n_{(0)})\\ &=&n_{(0)(0)}\otimes \beta^{-1}(n_{(0)(1)})\cdot
(beta^{-1}(S(n_{(1)}))\cdot m)\\ &=&n_{(0)(0)}\otimes \beta^{-1}(n_{(0)(1)}S(n_{(1)}))\cdot m\\ &=&n_{(0)}\otimes \beta^{-1}(n_{(1)(1)}S(n_{(1)(2)}))\cdot
m\\ &=&n_{(0)}\otimes \varepsilon(n_{(1)})m\\ &=& n\otimes m. \end{eqnarray*}

The fact that $ c_{M,N}^{-1}\circ c_{M,N}=id$ is similar.
 This completes the proof. \hfill $\blacksquare$
\\

 Let $H$ be a Hopf quasigroup and
 $G={\sl Aut}_{HQG}(H)
 \times {\sl Aut}_{HQG}(H)$.
 Define $\mathcal {YDQ}(H)$ as the
  disjoint union of all $_{H}\mathcal {YDQ}^{H}(\alpha,\beta)$
 with $(\alpha,\beta)\in G$. If we endow $\mathcal {YDQ}(H)$
 with tensor product shown in Proposition 3.1,
 then $\mathcal {YDQ}(H)$ becomes
 a monoidal category.

 Define a group homomorphism
 $\,\,\varphi: G\rightarrow Aut(\mathcal {YDQ}(H)),\quad
 (\alpha, \beta) \,\,\mapsto \,\,\varphi(\alpha,\beta)\,\,$
 on components as follows:
\begin{eqnarray*} \varphi_{(\alpha,\beta)}: {_{H}}\mathcal {YDQ}^{H}(\gamma,\delta)&\rightarrow& {_{H}}\mathcal {YDQ}^{H}((\alpha,\beta)\ast
(\gamma,\delta)\ast (\alpha,\beta)^{-1}),\\ \quad \quad \quad \quad \quad
 \quad \varphi_{(\alpha,\beta)}(N)&=& {}^{(\alpha,\beta)} N,
\end{eqnarray*} and the functor $\varphi_{(\alpha,\beta)}$ acts as identity on morphisms.\\

 The braiding in $\mathcal {YDQ}(H)$
 is given by the family $\{c_{M,N}\}$
 in Proposition 3.4.
 So we get the following main theorem of this article.
\\

{\bf Theorem 3.6.} $\mathcal {YDQ}(H)$ is a braided $T$-category over $G$. \\

\section*{ACKNOWLEDGEMENTS}
 The work was partially supported by the NSF of China (NO. 11371088 and NO.11571173), and the Fundamental Research Funds for the Central Universities  (NO.
 CXLX12-0067).

 \noindent Wei Wang,\\
Address: Department of Mathematics, Southeast University, Jiangsu Nanjing 210096, P. R. CHINA, Email: weiwang2012spring@yahoo.com.

 \noindent Shuanhong Wang,\\
Address: Department of Mathematics, Southeast University, Jiangsu Nanjing 210096, P. R. CHINA, Email: shuanhwang2002@yahoo.com.


\begin{thebibliography}{100}

  \bibitem{AFGS15} J. N. Alonso $\acute{A}$lvarez, J. M. Fern$\acute{a}$ndez Vilaboa, R. Gonz$\acute{a}$lez Rodr$\acute{1}$guez and C. Soneira Calvo.
      Projections and Yetter-Drinfel'd modules over
  Hopf (co)quasigroups, {\sl J.Algebra} 443, 153-199, 2015.

  \bibitem{BJ12} T. Brzezi$\acute{n}$ski, Z. Jiao, Actions of Hopf quasigroups, {\sl Comm. Algebra} 40,681-696, 2012.


  \bibitem{PD89} P. J. Freyd, D. N. Yetter. Braided compact closed categories with applications to low-dimensional topology. {\sl Adv. Math.} 77, 156-182,
      1989.


  \bibitem{KM10} J. Klim, S. Majid, Hopf quasigroups and the algebraic 7-sphere, {\sl J. Algebra} 323, 3067-3110, 2010.

  \bibitem{P07} J.M. P¨¦rez-Izquierdo, Algebras, hyperalgebras, nonassociative bialgebras and loops, {\sl Adv. Math.} 208, 834-876, 2007.

  \bibitem{PM07} P. Panaite, M. D. Staic. Generalized (anti) Yetter-Drinfel'd modules as components of a braided T-category. {\sl Israel J. Math.} 158,
      349-366, 2007.

  \bibitem{M69} M. E. Sweedler. Hopf Algebras. {\sl Benjamin, New York}, 1969.

  \bibitem{V08} V. G. Turaev. Crossed group-categories. {\sl Arab. J. Sci. Eng. Sect. C Theme Issues} 33(2C), 483-503, 2008.

  \bibitem{V94} V. G. Turaev. Quantum Invariants of Knots and $3$-Manifolds. {\sl de Gruyter Stud. Math.} de Gruyter, Berlin 18, 1994.


  \bibitem{M04} M. Zunino. Yetter-Drinfeld modules for crossed structures. {\sl J. Pure Appl. Algebra} 193, 313-343, 2004.

\end{thebibliography}
\end{document}